\documentclass[12pt]{amsart}
\usepackage[latin1]{inputenc}
\usepackage{mathptmx}
\usepackage{newalg}

\usepackage{pstcol}
\def\vertex{\pscircle[fillstyle=solid,fillcolor=black]{0.07}}
\definecolor{light}{gray}{0.9}
\definecolor{medium}{gray}{0.8}

\newtheorem{theorem}{Theorem}
\newtheorem{lemma}[theorem]{Lemma}

\theoremstyle{definition}
\newtheorem{remark}[theorem]{Remark}

\def\ZZ{{\mathbb Z}}
\def\RR{{\mathbb R}}

\def\cR{{\mathcal R}}

\def\inte{\operatorname{int}}
\def\Hilb{\operatorname{Hilb}}
\def\Supp{\operatorname{Supp}}
\def\gp{\operatorname{gp}}
\def\GL{\operatorname{GL}}
\def\id{\operatorname{id}}

\def\spec#1{\textsc{\texttt{#1}}}

\let\epsilon=\varepsilon

\begin{document}
\title{On the integral Caratheodory property}

\author{Winfried Bruns}
\address{Universit\"at Osnabr\"uck, FB Mathematik/Informatik, 49069
Osnabr\"uck, Germany} \email{wbruns@uos.de}

\begin{abstract}
In this note we document the existence of a finitely generated
rational cone that is not covered by its unimodular Hilbert
subcones, but satisfies the integral Carath\'eo\-do\-ry property. We
explain the algorithms that decide these properties and describe our
experimental approach that led to the discovery of the examples.
\end{abstract}

\maketitle

\section{Introduction}

Let $C\subset\RR^d$ be a finitely generated rational cone, i.~e.\
the set of all linear combinations $a_1x_1+\dots+a_nx_n$ of rational
vectors $x_1,\dots,x_d$ with coefficients from $\RR_+$. We can of
course assume that $x_i\in\ZZ^d$, $i=1,\dots,n$. In this note a
\emph{cone} is always supposed to be rational and finitely
generated. Moreover, we will assume that $C$ is \emph{pointed}: if
$x,-x\in C$, then $x=0$. Finally, it is tacitly understood that $C$
has full dimension $d$.

The monoid $M(C)=C\cap\ZZ^d$ is finitely generated by Gordan's lemma
(for example, see \cite[Section 2.A]{BGu:PRK}). Since $C$ is
pointed, $M$ is a \emph{positive} monoid so that $0$ is the only
invertible element in $M(C)$.

It is not hard to see that $M(C)$ has a unique minimal system of
generators that we call its \emph{Hilbert basis}, denoted by
$\Hilb(M(C))$ or simply $\Hilb(C)$. It consists of those elements
$z\neq 0$ of $M(C)$ that have no decomposition $z=x+y$ in $M(C)$
with $y,z\neq 0$.

We want to discuss combinatorial conditions on $\Hilb(C)$ expressing
that $C$ or $M(C)$ is covered by certain ``simple'' subcones or
submonoids, respectively. To this end we define a \emph{$u$-subcone}
of $C$ to be a subcone generated by vectors
$x_1,\dots,x_d\in\Hilb(C)$ that form a basis of the group $\ZZ^d$.
In particular, $x_1,\dots,x_d$ are linearly independent, and if just
this weaker condition is satisfied, then the cone $S$ generated by
$x_1,\dots,x_d$ is called an \emph{$f$-subcone}. In this case we let
$\Gamma(S)$ denote the subgroup of $\ZZ^d$ generated by
$x_1,\dots,x_d$ and $\Sigma(S)$ the submonoid of $\ZZ^d$ generated
by $x_1,\dots,x_d$. Note that $S$ is a $u$-subcone if and only if
$\Gamma(S)=\ZZ^d$, or, equivalently, $\Sigma(S)=S\cap\ZZ^d$.

One says that $C$ satisfies (UHC) if $C$ is the union of its
$u$-subcones. The letter U stands for \emph{unimodular}, H reminds
us of the condition that the generators of the $u$-subcones belong
to $\Hilb(C)$, and C simply stands for \emph{cover}.

A weaker condition than (UHC) is the \emph{integral Carath\'eodory
property} (ICP). One says that $C$ has (ICP) if every element of
$M(C)$ can be written as a linear combination of at most $d$
elements $x_i\in\Hilb(C)$ with \emph{integral} nonnegative
coefficients $a_i$. The terminology is motivated by Carath\'eodory's
theorem: let $y_1,\dots,y_m$ be a minimal system of generators of
the cone $C$; then every element $y\in C$ is a linear combination
$y=a_1y_{i_1}+\dots+a_dy_{i_d}$  with nonnegative real coefficients.

Both (UHC) and (ICP) can be formulated more generally for positive
affine monoids $M\subset\ZZ^d$. However, it is easy to see that
every monoid $M$ satisfying (UHC) is given in the form
$M=\RR_+M\cap\ZZ^d$. By a theorem of Bruns and Gubeladze
\cite[Theorem 6.1]{BGu:cov} the same holds true if $M$ satisfies
(ICP), provided the group $\gp(M)$ generated by $M$ equals $\ZZ^d$.
In loc.\ cit.\ it is also shown that (ICP) is equivalent to the
formally stronger condition that $M$ is the union of its submonoids
$\Sigma(S)$. (This condition is called (FHC) in \cite{BGu:cov}.) The
equivalence is crucial for our note, and therefore we reproduce the
statement and its proof in Theorem \ref{ICP=>FHC}.

While we view (UHC) and (ICP) as structural properties of (normal)
affine monoids, these properties have first been discussed in the
context of integer programming: see Cook, Fonlupt and Schrijver
\cite{CFS} and Seb\H{o} \cite{Sebo}.

It was asked by Seb\H{o} \cite{Sebo} whether every cone $C$ has
(ICP) or (UHC), and he proved that (UHC) holds if $d\le 3$. He
actually proved a stronger statement: $C$ has a triangulation by
$u$-subcones. A counterexample to (UHC) in dimension $6$, called
$C_{10}$ in the following, was found by Bruns and Gubeladze
\cite{BGu:cov}, and then verified to violate (ICP), too, in
cooperation with Henk, Martin, and Weismantel \cite{BGHMW}. Despite
the existence of the counterexample, one can fairly say, at least
heuristically, that almost all cones satisfy (UHC).

It remained an open problem whether (UHC) is strictly stronger than
(ICP). In this note we want to document the existence of cones that
satisfy (ICP) but fail (UHC), explain the algorithms that decide
(UHC) and (ICP), and describe the experimental approach that led to
the discovery of the examples.

All our experiments seem to indicate that $C_{10}$ is the core
counterexample to (ICP) and (UHC). In fact, all counterexamples to
these properties that we have been found contain it. It would be
very desirable indeed to clarify the situation in dimensions $4$ and
$5$.
\medskip

\emph{Acknowledgement.}\enspace The author is very grateful to
Joseph Gubeladze for inspiring discussions and to the Mathematisches
Forschungsinstitut Oberwolfach where the first steps of this project
were taken during a joint visit within the MFO's RiP program.

\section{Deciding UHC}

Let us say that $x\in C$ is \emph{$u$-covered} if it is contained in
a $u$-subcone. A subset of $C$ is $u$-covered if each of its
elements is $u$-covered. Using this simple terminology, we can
describe an algorithm deciding (UHC); see Table \ref{UHCalgo}.

\begin{table}[hbt]
\begin{algorithm}{\spec{unicover}}{D,n}
\begin{FOR}{i \= n \TO N}\\
\begin{IF}{D\subset U_i}
\RETURN
\end{IF}\\
\begin{IF}{\inte(D)\cap\inte(U_i)\neq\emptyset}
(D_1,D_2)\= \spec{split}(D,U_i)\\
\hspace*{0.07cm}\spec{unicover}(D_1,i)\\
\hspace*{0.07cm}\spec{unicover}(D_2,i) \\
\RETURN
\end{IF}
\end{FOR} \\
\spec{output}(\text{$D$ not $u$-covered})\\
\RETURN
\end{algorithm}
\begin{algorithm}{\spec{main}}{}
\text{Create the list $U_1,\dots,U_N$ of $u$-subcones of $C$}\\
\spec{unicover}(C,1)
\end{algorithm}
\caption{An algorithm deciding UHC}\label{UHCalgo}
\end{table}

In the algorithm we use a function named \spec{split}. It decomposes
$D$ along a support hyperplane $H$ of a $u$-subcone $U$ such that
$D\cap H^>\neq\emptyset$ as well as $D\cap H^<\neq\emptyset$. Such a
hyperplane does indeed exist if $\inte(D)\cap\inte(U)\neq\emptyset$,
but $D\not\subset U$. The cones produced are $D_1=D\cap H^+$ and
$D_2=D\cap H^-$. (The open halfspaces determined by $H$ are denoted
by $H^>$ and $H^<$, and $H^+$ and $H^-$ are the corresponding closed
halfspaces.)
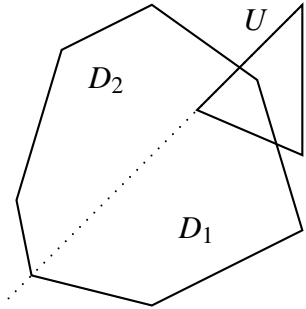
\begin{figure}[hbt]
\begin{center}
\psset{unit=2cm}
\begin{pspicture}(-1,-1)(1,1)
 \pspolygon(-0.8,-0.8)(0,-1)(1,-0.5)(0.7,0.5)(0,1)(-0.6,0.7)(-0.9,-0.3)
 \pspolygon(0.3,0.3)(1,1)(1.0,0)
 \psline[linestyle=dotted](1,1)(-1,-1)
 \rput(-0.3,0.5){$D_2$}
 \rput(0.3,-0.5){$D_1$}
 \rput(0.7,0.9){$U$}
\end{pspicture}
\end{center}
\caption{The function \spec{split}} \label{Figsplit}
\end{figure}

It is easy to see that the algorithm terminates: there are only
finitely many hyperplanes by which we split subcones. Therefore only
finitely many subcones can be created.

In order to check the correctness of the algorithm, observe that at
each step in the \textbf{for} loop in \spec{unicover} none of the
$u$-subcones $U_1,\dots,U_{i-1}$ of $C$ intersects the interior of
$D$. In fact, the index $i$ is only increased if $\inte(D)\cap
\inte(U_i)=\emptyset$. (In the recursive call the index $i$ that has
been reached at the parent level starts the \textbf{for} loop at the
child level.)

Thus, when the loop terminates with $i=N+1$, none of the
$u$-subcones of $C$ intersects the interior of $D$. So $D$, and
consequently $C$, indeed contains a vector that is not $u$-covered.
Conversely, if $C$ contains such a vector $x$, then on each level of
the recursion tree one finds a subcone $D$ containing $x$, and for
such $D$ the condition $D\subset U_i$ can never be satisfied.
Therefore there exists an end node of the recursion tree at which
the loop is left with $i=N+1$.

The pseudocode in Table \ref{UHCalgo} is a somewhat simplistic
sketch of the actual implementation since it is necessary to cope
with substantial memory requirements. For example, the list
$U_1,\dots,U_N$ is not produced a priori, but extended whenever
necessary and always kept as small as possible. Moreover, all
allocated memory is recycled carefully within the program.

Instead of starting the covering algorithm with the full cone $C$,
the actual implementation uses the output of a preprocessor that
computes several triangulations $\Delta_1,\dots,\Delta_t$ of $C$,
The input to \spec{unicover} (and also to \spec{caradec} below) is
the list of intersections $D_1\cap\dots\cap D_t$ where $D_i$ is a
nonunimodular simplicial cone in $\Delta_i$.

\section{Deciding ICP}\label{ICP}

Let us first fix some terminology that parallels that for (UHC). An
element $x\in C\cap \ZZ^d$ is \emph{$f$-covered} if it belongs to
one of the monoids $\Sigma(S)$ where $S$ is an $f$-subcone, and a
subset of $C$ is $f$-covered if each of its elements is $f$-covered.

The lemma contains the basic criterion by which we can check that
$C$ is $f$-covered. In the theorem following it, we will then see
that this property is equivalent to (ICP).

\begin{lemma}\label{ICPcrit}\leavevmode
\begin{itemize}
\item[(a)] Let $G_1,\dots,G_n$ be subgroups of $\ZZ^d$, and let $N$
be a residue class of $\ZZ^d$ modulo $G_1\cap\dots\cap G_n$. Then
$N\subset G_1\cup\dots\cup G_n$ if and only if $N\cap
(G_1\cup\dots\cup G_n)\neq \emptyset$.

\item[(b)] Let $S_1,\dots,S_n$ be $f$-subcones of $C$, each containing
the $d$-dimensional subcone $D$ of $C$. If every residue class of
$\ZZ^d$ modulo $\Gamma(S_1)\cap\dots\cap\Gamma(S_n)$ meets
$\Gamma(S_1)\cup\dots\cup\Gamma(S_n)$, then $D$ is $f$-covered.

\item[(c)] Let $D$ be a $d$-dimensional subcone of $C$ with the
following property: for every $f$-subcone $S$ either $D\subset S$ or
$\inte(D)\cap\inte(S)=\emptyset$. Furthermore let $G_D$ be the
intersection of the groups $\Gamma(S)$, $S\supset D$, and $H_D$
their union. Then $D$ is $f$-covered if and only if every residue
class of $\ZZ^d$ modulo $G_D$ meets $H_D$.
\end{itemize}
\end{lemma}

\begin{proof}
(a) Suppose that $N\cap(G_1\cup\dots\cup G_n)\neq \emptyset$, and
let $x$ be an element in the intersection, $x\in G_i$. The subgroup
$G'=G_1\cap\dots\cap G_n$ is contained in $G_i$, and so
$N=x+G'\subset G_i$. The converse implication is trivial.

(b) Let $x\in D\cap\ZZ^d$. It follows from (a) that
$x\in\Gamma(S_i)$ for some $i$. But $x\in S_i$, too. Therefore
$x\in\Gamma(S)\cap S_i=\Sigma(S_i)$. (At this point we use that the
generators of $S_i$ are linearly independent.)

(c) It only remains to show the necessity of the condition. For it
we only need to observe that every residue class $N$ of $\ZZ^d$
modulo $G_D$ meets $\inte(D)$. By hypothesis on $D$, an element
$x\in N\cap\inte(D)$ is $f$-covered if and only if $x\in \Sigma(S)$
for some $f$-subcone $S$ containing $D$.
\end{proof}

We include the next theorem and its proof for the convenience of the
reader. It is a simplified version of \cite[Theorem 6.1]{BGu:cov}
whose proof contains the crucial ideas for the algorithm deciding
(ICP).

\begin{theorem}\label{ICP=>FHC}
Let $M\subset \ZZ^d$ be a positive affine monoid such that
$\gp(M)=\ZZ^d$. If $M$ satisfies \emph{(ICP)}, then
$M=\RR_+M\cap\ZZ^d$, and every element of $M$ is $f$-covered.
\end{theorem}

\begin{proof}
We dissect $C=\RR_+M$ along all the support hyperplanes of the cones
spanned by linearly independent vectors $x_1,\dots,x_d$ of
$\Hilb(M)$ into \emph{elementary} subcones. Set $\bar
M=\RR_+M\cap\ZZ^d$ and choose $x\in \bar M$. Suppose that $x$ has no
representation as a linear combination $x=a_1y_1+\dots+a_dy_d$,
$a_1,\dots,a_d\in\ZZ_+$ and $y_1,\dots,y_d\in\Hilb(M)$ linearly
independent. The element $x$ belongs to one of the elementary
subcones $D$, and as in the proof of the lemma it follows that there
exists a finite index subgroup $G$ of $\ZZ^d$ such that no element
$z$ of $(x+G)\cap\inte(D)$ has a representation
$a_1y_1+\dots+a_dy_d$ with $a_1,\dots,a_d\in\ZZ_+$ and
$y_1,\dots,y_d\in\Hilb(M)$ linearly independent.

The crucial point is that $\bar M\setminus M$ is contained in the
union of finitely many hyperplanes (see \cite[Section
2.B]{BGu:PRK}), and the same applies to all elements of $M$ that are
linear combinations of linearly dependent elements of $\Hilb(M)$.
But $(x+G)\cap\inte(D)$ is not contained in the union of finitely
many hyperplanes, and so must contain elements of $M$. This is
impossible if $M$ satisfies (ICP).
\end{proof}

\begin{table}[hbt]
\begin{algorithm}{\spec{caradec}}{D,G,\cR,n}
\begin{FOR}{i\= n \TO N}\\
\begin{IF}{D\subset S_i}
\cR'\=\emptyset\\
\begin{FOR}{\bigl(x\in \cR,\ y\in G/(G\cap\Gamma(S_i))\bigr)}\\
\begin{IF}{x+y\notin\Gamma(S_i)}
\cR'=\cR'\cup\{x+y\}
\end{IF}
\end{FOR}\\
G\= G\cap \Gamma(S_i),\quad
\cR\=\cR'\\
\begin{IF}{\cR=\emptyset}
\RETURN
\end{IF}
\end{IF}\\
\begin{IF}{D\not\subset S_i\ \algkey{and}\inte(D)\cap\inte(S_i)\neq\emptyset}
(D_1,D_2)\=\spec{split}(D,S_i)\\
\hspace*{0.07cm}\spec{caradec}(D_1, G ,\cR, i)\\
\hspace*{0.07cm}\spec{caradec}(D_2, G ,\cR, i)\\
\RETURN
\end{IF}
\end{FOR}\\
\spec{output}(\text{$D$ not $f$-covered})\\
\RETURN
\end{algorithm}
\begin{algorithm}{\spec{main}}{}
\text{Create the list $S_1,\dots,S_N$ of $f$-subcones of $C$}\\
\spec{caradec}(C,\ZZ^d,\{0\},1)
\end{algorithm}
\caption{An algorithm deciding ICP}\label{ICPalgo}
\end{table}

For the algorithm deciding (ICP) we have to enrich our data
structure by those components that have shown up in the proof of the
lemma. Subcones are replaced by triples $(D,G,\mathcal R)$ where $D$
is a subcone of $C$, $G$ is a finite index subgroup of $\ZZ^d$, and
$\mathcal R$ is a list of residue classes in $\ZZ^d/G$. In $\cR$
each residue class is represented by a single vector that belongs to
it, and in the algorithm (see Table \ref{ICPalgo}) the loop
$$
\textbf{for}\ \bigl(x\in \cR,\ y\in G/(G\cap\Gamma(S_i))\bigr)
$$
runs over all elements of $\cR\times G/(G\cap\Gamma(S_i))$.

Again it is clear that the algorithm terminates after finitely many
steps: the number of hyperplanes that we can use to split subcones
of $C$ is still finite (though larger than for (UHC)).

The crucial point for {\tt caradec} is that at each step in the loop
the $f$-cones $S_1,\dots,S_N$ satisfy the following conditions:
\begin{enumerate}
\item for each $j\le i-1$ either $S_j\supset D$ or
$\inte(S_j)\cap\inte(D)=\emptyset$;

\item $G$ is the intersection of all groups $\Gamma(S_j)$, $j\le
i-1$, for which $D\subset S_j$;

\item $\cR$ is the list of those residue classes in $\ZZ^d/G$ that
are not contained in the union of the groups $\Gamma(S_j)$, $j\le
i-1$, $D\subset S_j$.
\end{enumerate}

We have only to check that these conditions remain satisfied when
$S_i$ is tested against $D$. To this end let $S_{k_1},\dots,S_{k_m}$
be those among $S_1,\dots,S_{i-1}$ that contain $D$.

If $\inte(S_i)\cap \inte(D)=\emptyset$, then $D\not\subset S_i$, and
this case is done.

If $\inte(S_i)\cap \inte(D)\neq\emptyset$, but $D\not\subset S_i$,
then $i$ is not increased (!) and all three conditions are inherited
by both $D_1$ and $D_2$: among the $S_j$, $j\le i-1$, exactly
$S_{k_1},\dots,S_{k_m}$ contain $D_1$ or $D_2$, simply because
$S_j\supset D_1$ or $S_j\supset D_1$ implies $\inte(S_j)\cap
\inte(D)\neq \emptyset$, and so $S_j\supset D$.

But if $S_i\supset D$, the bookkeeping is also correct. Evidently
$G$ is replaced by the correct group $G\cap\Gamma(S_i)$. Next
observe that all residue classes of $G\cap\Gamma(S_i)$ that are
contained in residue classes modulo $G$ not appearing in $\cR$
remain in $\Gamma(S_{k_1})\cup\dots\cup \Gamma(S_{k_m})$. On the
other hand, those that refine elements of $\cR$ must belong to
$\Gamma(S_i)$ to be in $\Gamma(S_{k_1})\cup\dots\cup
\Gamma(S_{k_m})\cup \Gamma(S_i)$. The correctness of the algorithm
follows now immediately from Lemma \ref{ICPcrit}.

The biggest hurdle for it are the lists $\cR$ of residue classes
that usually become extremely long already in dimension $6$.
Moreover, along each branch of the recursion tree, several of them
must be kept in memory. (This problem cannot be eliminated by a
nonrecursive implementation.)

The growth of the list $\cR$ can be estimated. Set
$$
e=\#\bigl(G/(G\cap\Gamma(S_i))\bigr)\quad\text{and}\quad
e'=\#(\ZZ^d/\Gamma(S_i)).
$$
Each element $x\in\cR$ is involved in $e$ vectors $x+y$. At most one
of them lies in $\Gamma(S_i)$ since the vectors $y$ belong to
pairwise different residue classes modulo $\Gamma(S_i)$. Therefore
$$
\#(\cR')\ge (e-1)\#(\cR).
$$
If the elements of $\cR$ are randomly distributed over the residue
classes of $\ZZ^d$ modulo $\Gamma(S_i)$, then the expected share of
vectors $x+y\in\Gamma(S_i)$ drops to $1/e'$.

Instead of keeping the lists of residue classes in memory, one could
alternatively try to follow the recursion tree along the whole list
$S_1,\dots,S_N$, compute only $G$ along each branch and test the
residue classes one by one only at the end nodes. However, this
approach seems unfeasible since it derives no advantage from the
case $e=1$, which fortunately happens frequently and often stops the
recursion before the end of $S_1,\dots,S_N$ is reached.

The list $S_1,\dots,S_N$ is actually scanned in growing order of the
determinants of the $S_i$. This has turned out very effective, at
least for those cones that satisfy (ICP). In fact, all cones in
Table \ref{TUC} with (ICP) are covered by $f$-subcones of
determinant $\le 2$.

In addition to \spec{caradec} we use a Monte Carlo approach for
disproving (ICP). It reads the output of \spec{unicover}, computes a
large number of vectors in the non-$u$-covered subcones of $C$ and
tests whether they are $f$-covered.

\begin{remark}\label{ratio}
\spec{caradec} provides us with a precise measure for the failure of
(ICP), namely the ratios $\#(\cR)/\#(G)$ at the end nodes of the
recursion tree. For the cone $C_{10}$ (Table \ref{HilbC10}) there is
precisely one end node with $\cR\neq\emptyset$, and the ratio is
$32/15552=1/486$. The number of non-$f$-covered vectors in the Monte
Carlo test confirms the ratio rather precisely.

For the cone $C_{15}''$ (Table \ref{TUC}) there is again a single
group with $\cR\neq\emptyset$ and ratio
$$
2,4468,480/2,286,144,000,000\approx 1.07/10^6.
$$
So the Monte Carlo test cannot be expected to be conclusive with
${<10^6}$ test vectors.
\end{remark}

\section{The search}\label{hunt}

Let us recapitulate an important notion from \cite{BGu:cov}. An
element $x$ of $\Hilb(C)$ is called \emph{destructive} if
$H'=\Hilb(C)\setminus\{x\}$ is not the Hilbert basis of $\RR_+H'$.
We say that $C$ is \emph{tight} if every element of $\Hilb(C)$ is
destructive. The crucial role of tight cones for (UHC) and (ICP) is
illuminated by the following lemma \cite[Corollary 2.3]{BGu:cov}.

\begin{lemma}\label{tight}
Let $C$ be a cone that is a counterexample to \emph{(UHC)} or
\emph{(ICP)}. Suppose $C$ is minimal first with respect to dimension
and second with respect to $\#\Hilb(C)$. Then $C$ is tight.
\end{lemma}

\begin{remark}
Updating the information in \cite{BGu:cov} we mention that tight
cones exist in all dimensions $d\ge 3$. The first $3$-dimensional
tight cone was found by P. Dueck. The smallest such cone found by
the author has a Hilbert basis of $19$ elements. The elements of the
Hilbert basis in the extreme rays form a regular hexagon (with
respect to the action of $\GL_3(\ZZ)$) so that the cone has the
dihedral group $D_6$ as its automorphism group. The regularity is in
indication that it may be the smallest possible tight cone. (Here
and in the following the automorphism group of a cone $C$ is always
understood to be the automorphism group of the monoid $C\cap\ZZ^d$.)
\end{remark}

Our search for counterexamples has been based on the crucial Lemma
\ref{tight}. We produce a set of random vectors, consider them as
the generating set of a cone $C$, and then use a program named
\spec{shrink} to remove nondestructive elements of $\Hilb(C)$ until
a tight cone is reached. (\spec{shrink} is based on the same
algorithm as \spec{normaliz}; see \cite{BK:nmz, BK:comp}.) Almost
always, $C$ shrinks to the $0$-cone, but sometimes a nontrivial
tight cone emerges. Then \spec{unicover}, and possibly
\spec{caradec}, are invoked.

When we started the search in spring 1998, we used cones over
randomly generated lattice parallelepipeds. In May 1998 the search
stopped with the counterexample $C_{10}$. Its Hilbert basis is shown
in Table \ref{HilbC10}. The cone $C_{10}$ has $27$ support
hyperplanes.
\begin{table}[hbt]
\begin{alignat*}{2}
 z_1 & = (0,\,1,\,0,\,0,\,0,\,0) & \qquad z_6 &=
(1,\,0,\,2,\,1,\,1,\,2) \\
 z_2 & = (0,\,0,\,1,\,0,\,0,\,0) & \qquad z_7 &=
(1,\,2,\,0,\,2,\,1,\,1) \\
 z_3 & = (0,\,0,\,0,\,1,\,0,\,0) & \qquad z_8 &=
(1,\,1,\,2,\,0,\,2,\,1) \\
 z_4 & = (0,\,0,\,0,\,0,\,1,\,0) & \qquad z_9 &=
(1,\,1,\,1,\,2,\,0,\,2) \\
 z_5 & = (0,\,0,\,0,\,0,\,0,\,1) & \qquad z_{10} &=
(1,\,2,\,1,\,1,\,2,\,0)
\end{alignat*}
\caption{$\Hilb(C_{10})$}\label{HilbC10}
\end{table}

The reader should note that for the questions considered in this
note we can always replace a given cone $C$ by $\phi(C)$ where
$\phi$ is an arbitrary transformation in $\GL_d(\ZZ)$. In this
sense, $C$ stands for a class of cones that are isomorphic under an
integral isomorphism of $\RR^d$. We express this fact by speaking of
different embeddings of a cone $C$.

While \spec{unicover} showed that $C_{10}$ fails (UHC), it was then
verified in cooperation with Henk, Martin, and Weismantel that
$C_{10}$ also fails (ICP). (\spec{caradec} was not written before
September 2006.) See Bruns and Gubeladze \cite{BGu:cov} and Bruns et
al.\ \cite{BGHMW} for more information on $C_{10}$.

The automorphism group of $C_{10}$ is remarkably large: it is the
Frobenius group $F_{20}$ of order 20, which acts transitively on
$z_1,\dots,z_{10}$. ($F_{20}$ is the semidirect product of $\ZZ_5$
with its automorphism group $\ZZ_5^*\cong \ZZ_4$.) From the
embedding above one can see that at least the dihedral group
$D_5\subset F_{20}$ is acting on $C_{10}$. All the remaining $10$
automorphisms have order $4$ and swap $\{z_1,\dots,z_5\}$ with
$\{z_6,\dots,z_{10}\}$. Moreover, $z_1,\dots,z_{10}$ all lie in the
hyperplane given by $-5\zeta_1+\zeta_2+\dots+\zeta_6=1$. The convex
hulls of $\{z_1,\dots,z_5\}$ and $\{z_6\dots,z_{10}\}$ are both
simplices of dimension $4$.

\begin{remark}
It was communicated to us by F. Santos that the lattice polytope
spanned by $\Hilb(C_{10})$ is a projection of the Ohsugi-Hibi
polytope \cite{OhH:nrt1}. The projection leads to the following
description of $C_{10}\cap\ZZ^6$. Consider the complete graph $K_5$
and decompose it into $2$ cycles of length $5$ as shown in Figure
\ref{K_5}.
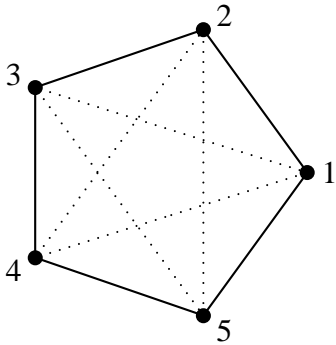
\begin{figure}[hbt]
\begin{center}
\psset{unit=2cm}
\def\vertex{\pscircle[fillstyle=solid,fillcolor=black]{0.05}}
\begin{pspicture}(-1,-1)(1,1)
 \def\A{1,0}
 \def\B{0.309,0.951}
 \def\C{-0.809,0.566}
 \def\D{-0.809,-0.566}
 \def\E{0.309,-0.951}
 \rput(\A)\vertex
 \rput(\B)\vertex
 \rput(\C)\vertex
 \rput(\D)\vertex
 \rput(\E)\vertex
 \pspolygon(\A)(\B)(\C)(\D)(\E)
 \pspolygon[linestyle=dotted](\A)(\C)(\E)(\B)(\D)
 \rput(1.15,0){$1$}
 \rput(0.45,1.05){$2$}
 \rput(-0.95,0.65){$3$}
 \rput(-0.95,-0.65){$4$}
 \rput(0.45,-1.05){$5$}
\end{pspicture}
\end{center}
\caption{Cycle decomposition of $K_5$} \label{K_5}
\end{figure}
Now choose the incidence vectors $(1,1,0,0,0)$ etc.\ of the edges in
the first cycle and prefix them with $0$. Then prefix the incidence
vectors $(1,0,1,0,0)$ etc.\ of the second cycle with $1$. The
resulting $10$ vectors in $\ZZ^6$ generate a monoid $M$ isomorphic
with $C_{10}\cap\ZZ^6$. While this description is even more
aesthetic than the one in Table \ref{HilbC10}, it has the
disadvantage that $\gp(M)$ is of index $2$ in $\ZZ^6$.
\end{remark}

In the summer of 1998 a second counterexample $C_{12}$ to (UHC) and
(ICP) emerged. It has a Hilbert basis of 12 elements. We continued
the search for two more years. The frustrating outcome was that
$C_{10}$ appeared over and over again, but no new counterexample
showed up (and even $C_{12}$ did not return until November 22,
2006).

The project was taken up again at the end of 2004 when our
department had installed a dual processor Opteron system with very
fast integer arithmetic. Nevertheless, the outcome of the search
remained as disappointing as it had been before.

Finally, in August 2006 we did what should have been done long
before, namely compare $C_{12}$ with $C_{10}$: it turned out that
$\Hilb(C_{12})$ (in an embedding that had to be found!) extends
$\Hilb(C_{10})$ by two vectors. Relative to $C_{10}$, this finding
explained why $C_{12}$ fails (UHC) and (ICP), too: the extra
$u$-subcones and $f$-subcones are not sufficient to cover all
integral vectors in $C_{10}$. (It also shows that one cannot speed
up shrinking by removing two vectors at a time.)

However, this not very surprising \emph{a posteriori} insight made
it suddenly clear that there might be many interesting objects in
the vicinity of $C_{10}$. Especially, when we approach $C_{10}$
along a shrink path, why should the stronger property (UHC) not be
lost before (ICP)? After a modification of \spec{shrink} we also
applied \spec{unicover} to the, say, $6$ last non-tight
approximations to $C_{10}$, and within hours many new non-(UHC)
cones emerged. Several of them defeated all Monte Carlo attacks on
(ICP). It became clear that \spec{caradec} had to be implemented,
and it indeed recognized many non-(UHC), but (ICP) cones.

Ironically, within a few weeks after we had given up our
narrow-minded insistence on checking only tight cones, two new
non-(UHC) such cones surfaced, both of them satisfying (ICP). They
appear as $C_{12}'$ and $C_{15}$ in Table \ref{moreHilb}.
\begin{table}[hbt]
\begin{alignat*}{4}
C_{12}:\quad&& z_{11}' & = (2,\,2,\,1,\,4,\,1,\,3) & \qquad
C_{15}:\quad&&
w_1 &=(2,\,1,\,0,\,5,\,1,\,5) \\
&& z_{12}' & = (2,\,3,\,1,\,4,\,1,\,2) & \qquad &&
w_2 &=(1,\,0,\,-1,\,4,\,0,\,4 ) \\
&&  && \qquad &&
w_3 &=(0,\,0,\,-1,\,1,\,0,\,1) \\
C_{12}':\quad&& z_{11}'' & = (0,\,-1,\,2,\,-1,\,-1,\,2) & \qquad &&
w_4&=(2,\,1,\,2,\,3,\,2,\,4) \\
&& z_{12}'' & = (1,\,0,\,3,\,0,\,0 ,\,3) & \qquad && w_5
&=(1,\,1,\,0,\,3,\,1,\,2).
\end{alignat*}
\caption{Additional vectors in $\Hilb(C_{12})$, $\Hilb(C_{12}')$,
$\Hilb(C_{15})$}\label{moreHilb}
\end{table}
Since all these cones contain $C_{10}$, we list only the extra
vectors that complement the Hilbert basis of $C_{10}$ (using the
embedding given in Table \ref{HilbC10}). The numbers of support
hyperplanes are $39$ for $C_{12}$, $40$ for $C_{12}'$, and $36$ for
$C_{15}$.

The most interesting after $C_{10}$ undoubtedly is $C_{12}'$, not
only because it satisfies (ICP). Its automorphism group -- certainly
invisible from the embedding given -- is again the Frobenius group
$F_{20}$. It is clear that $F_{20}$ cannot act transitively on
$\Hilb(C_{12}')$, which rather decomposes into an orbit of $10$
elements and one of $2$. However, this is by no means an extension
of the action of $F_{20}$ on $C_{10}$ since the orbit of two
elements is $\{z_1,z_5\}$! Only a subgroup isomorphic to
$\ZZ_2\times\ZZ_2$ restricts to $C_{10}$, showing that there are $5$
conjugate embeddings of $C_{10}$ into $C_{12}'$,  and each of them
contains $z_1,z_{5}$.

The convex hulls of $\{z_1,\dots,z_5,z_{11}''\}$ and
$\{z_6\dots,z_{10},z_{12}''\}$ are both bipyramids over a
tetrahedron. The bipyramids are situated in the parallel hyperplanes
with the equations $\zeta_1=0$ and $\zeta_1=1$. Both $C_{12}$ and
$C_{12}'$ have their Hilbert bases in the hyperplane spanned by
$\Hilb(C_{10})$ but this is not true for $C_{15}$.

Table \ref{TUC} lists all the $11$ tight non-(UHC) cones of
dimension $6$ that have been found by December 17, 2006, including
those mentioned already. In the last column we indicate whether the
Hilbert basis is contained in a hyperplane.
\begin{table}[hbt]
\def\strut{\vphantom{\Large I}}
\begin{tabular}{|c|c|c|c|c|c|}
\hline \strut &$\#\Hilb$&$\#\Supp$&ICP&Aut&flat\\\hline
\strut$C_{10}$&10&27&no&$F_{20}$&yes\\\hline
\strut$C_{12}$&12&39&no&$\ZZ_2\times\ZZ_2$&yes\\\hline
\strut$C_{12}'$&12&40&yes&$F_{20}$&yes\\\hline
\strut$C_{14}$&14&34&yes&$\{\id\}$&yes\\\hline
\strut$C_{14}'$&14&39&no&$\{\id\}$&yes\\\hline
\strut$C_{14}''$&14&42&yes&$\ZZ_2$&no\\\hline
\strut$C_{15}$&15&36&yes&$\{\id\}$&no\\\hline
\strut$C_{15'}$&15&36&yes&$\{\id\}$&yes\\\hline
\strut$C_{15}''$&15&44&no&$\{\id\}$&yes\\\hline
\strut$C_{16}$&16&49&no&$\ZZ_2$&yes\\\hline
\strut$C_{16}'$&16&36&no&$\ZZ_2$&yes\\\hline
\end{tabular}
\bigskip \caption{Tight non-(UHC) cones}\label{TUC}
\end{table}

The smallest non-(UHC), but (ICP) cone we have found has a Hilbert
basis of 11 elements, and the largest has a Hilbert basis of 24
elements (most likely (ICP)). Like all the others they are
extensions of $C_{10}$ which seems to be the core obstruction to
(ICP) and (UHC).

While we have the implications (UHC) $\implies$ (ICP) $\implies$
$M=\RR_+M\cap \ZZ^d$ (under the condition that $\gp(M)=\ZZ^d$), it
is now clear that the converse implications do not hold. However, it
remains an open problem whether all cones in dimensions $4$ and $5$
have (ICP) or even (UHC).

\section{Computational issues}\label{Comp}
All programs have been written in C. The tight cones in Section
\ref{hunt} were found by the Opteron (O) system mentioned above. It
runs Linux, and the executables have been produced by the gcc
compiler. In the following we will also mention computations on two
other systems, the author's Intel Core2 6600 (C2) with Windows XP
and the DJGPP port of gcc, and the University of Osnabrück's Itanium
(I) system with Linux and the Intel compiler icc. The machines (O)
and (C2) are close to each other in speed; (I) is somewhat slower
(for integer arithmetic), but has very large memory (32 GB).

Some of the equipment used in the 1998 computations is still
accessible. This allowed us to measure the gain in speed by improved
hardware: the factor is $\ge 40$. Moreover, a better implementation
of \spec{shrink} yields an acceleration by a factor $\ge 3$. In
other words, $4$ months of the 1998 search take now a single day.

The number of cones shrunk by \spec{shrink} per second depends very
much on the parameters used for their creation. The performance for
$6$-dimensional cones generated by random $0$-$1$-vectors, whose
number varies between $6$ and $26$, is about $1000$ per second on
(O) or (C2). Cones over $5$-dimensional parallelotopes of Euclidean
volume $\le 30$ are shrunk at a rate of $0.6$ per second. The output
of tight cones is nevertheless comparable.

The cones in Section \ref{hunt} are light food for \spec{unicover}.
For example, the running time for $C_{14}''$ on (C2) is $1.7$
seconds. About $2.5$ million vectors are created, but the list of
vectors in memory simultaneously is bounded by $15,000$.

While all the other programs use 32 bit arithmetic, \spec{caradec}
is set to 64 bit. It has no problem with all the cones mentioned, as
long as they have (ICP), simply because in all cases they are
$f$-covered by cones of determinant $\le 2$. The running time for
$C_{14}''$ on (C2) is $6.9$ seconds.

Despite of its sometimes enormous appetite for memory,
\spec{caradec} has also been successful for all the cones of Table
\ref{TUC} that lack (ICP), with the exception of $C_{16}$ and
$C_{16}'$. It failed for these cones though it was allowed $200$
million vectors in memory.

The (ICP) property of $C_{16}$ and $C_{16}'$ was falsified by the
Monte Carlo method with $1$ million test vectors for each of the
non-(UHC) subcones produced by \spec{unicover} ($17.4$ seconds on
(C2) for $C_{16}$).

The longest successful run of \spec{caradec} with a negative result
was $C_{15}''$: $1,990$ seconds on (I), $2.1$ billion vectors, $110$
million simultaneously. The Monte Carlo method does the job with $1$
million vectors for the single non-(UHC) subcone in $2.7$ seconds on
(C2). See Section \ref{ICP} for further data $C_{15}''$.

\end{document}